# Unified approach to discretization of flow in fractured porous media

J. M. Nordbotten[1,2], W. M. Boon[1], A. Fumagalli[1], E. Keilegavlen[1]

[1] Department of Mathematics, University of Bergen, N-5020 Bergen, Norway
[2] Princeton Environmental Institute, Princeton University, Princeton, NJ 08544, USA

## Abstract

In this paper, we introduce a mortar-based approach to discretizing flow in fractured porous media, which we term the mixed-dimensional flux coupling scheme. Our formulation is agnostic to the discretizations used to discretize the fluid flow equations in the porous medium and in the fractures, and as such it represents a unified approach to integrated fractured geometries into any existing discretization framework. In particular, several existing discretization approaches for fractured porous media can be seen as special instances of the approach proposed herein.

We provide an abstract stability theory for our approach, which provides explicit guidance into the grids used to discretize the fractures and the porous medium, as dependent on discretization methods chosen for the respective domains. The theoretical results are sustained by numerical examples, wherein we utilize our framework to simulate flow in 2D and 3D fractured media using control volume methods (both two-point and multi-point flux), Lagrangian finite element methods, mixed finite element methods, and virtual element methods. As expected, regardless of the ambient methods chosen, our approach leads to stable and convergent discretizations for the fractured problems considered, within the limits of the discretization schemes.

## 1. Introduction

Flow in porous media with thin inclusions is an important process both within subsurface and industrial materials. Our main focus herein is on the subsurface, where the thin inclusions represent fractures, and the fracture space can be either open or filled. We will thus simply refer to fractured porous media in what follows. However, thin inclusions may also be engineered in artificial porous media for the purpose of fluid flow control.

Fluid flow in fractured porous media is a dominating process in several subsurface applications, ranging from geothermal energy production, shale gas recovery and nuclear waste deposits. As such, accurate and reliable numerical representations have been an important topic of research, and Rainer Helmig has been a major contributor to the field for more than three decades. Existing discrete representations of fractured porous media fall in two major categories, depending on whether the fractures conform to the underlying discrete grid representing the porous materials. So-called "unfitted" discretizations, wherein the fractures are allowed to be arbitrary with respect to the grid, have seen significant research and developments in recent years (see e.g. [1, 2]). Our focus herein is, in contrast, on discretizations where the discrete grid resolves the fractures, which are conceptually simpler than unfitted discretizations.

Early research into numerical simulation and conforming discretization of fractured porous media was spear-headed by among others Rainer Helmig and his collaborators [3]. This early work was centered around lowest-order finite element discretizations. Later, it was understood that local conservation properties were important for discretization methods for flow in porous media, and conforming discretizations of fractured porous media were developed based on control volume approaches [4, 5], mixed finite element methods [6, 7], mimetic finite differences [8] and virtual element methods [9]. See also [10] for a comparison study.

A recent development in the mathematical representation of fractured porous media is the modeling and interpretation of fractures as lower-dimensional manifolds [11, 7, 12]. This concept allows for the introduction of mixed-dimensional partial differential equations (md-PDEs), wherein partial differential equations are defined, in a coupled sense, both in the porous material, lower-dimensional fractures, and yet lower-dimensional intersections. In this abstraction, it can be shown that the mathematical models for fractured porous media, can be cast in a rich functional-analysis framework, ensuring well-posedness, and thus existence and uniqueness, of solutions [13].

In this manuscript, we revisit conforming discretizations of fractured porous media within the context of md-PDEs. We show, by introducing explicit coupling variables in the spirit of mortar methods [14, 11, 7, 15], an abstract framework for constructing a conforming fracture discretization from any discretization of non-fractured porous media. We term this approach the mixed-dimensional flux-coupling (MDFC) method. Viewed from the discretization within each dimension, the coupling between dimensions takes the form of standard boundary value problems, thus any implementation that can account for Dirichlet and Neumann boundary data can be applied to fractured media with minimal adaptations. Our approach thus unifies the various previous developments reviewed above.

We concretize the abstract framework by applying it to well-known discretizations from literature, establishing (in some cases for the first time) that these discretizations are well-posed. To illustrate the versatility of the framework, we provide numerical examples showing how five different discretization methods for non-fractured porous media can be applied as discretization methods for fractured porous media. Of these discretizations, when using mixed finite elements or standard finite elements for the non-fractured media, we recover earlier methods referenced above. In the case of finite volume (both two-point and multi-point flux) and virtual element methods, our approach effectively leads to a discretization scheme not previously discussed in literature. Our numerical examples, which include a 2D case where we use non-matching grids between the dimensions and a relatively complex 3D case, highlight the convergence properties and stability of MDFC even for degenerating parameters.

The remaining manuscript honors the following structure: In section 2, we introduce our novel approach to unifying discretization methods for fractured media. Thereafter, in section 3, we show the stability of the approach theoretically, which emphasizes the conditions required between the (in principle non-matching) grids discretizing the matrix and fractures. Numerical examples and verification are presented before concluding the paper.

## 2. Modeling fractured porous media

In this section we introduce our model for fractured media, first by a single fracture, and then extended to general fracture networks.

### 2.1. Domain with a single fracture

Flow in (fractured) porous media can lead to complex and non-linear governing equations. However, at the heart usually lies a second-order partial differential equation, which upon linearization (i.e. within a Newton iteration) thus takes the classical form for a pressure $p^3$ and flux $q^3$

$$\nabla \cdot q^3 + \psi^3 = 0 \quad \text{on} \quad \Omega^3 \quad (2.1)$$

$$-\kappa \nabla p^3 = q^3 \quad \text{on} \quad \Omega^3 \quad (2.2)$$

$$q^3 \cdot n^3 = \lambda_+^2 \quad \text{on} \quad \partial_{\Omega_+^2} \Omega^3 \quad (2.3)$$

$$q^3 \cdot n^3 = \lambda_-^2 \quad \text{on} \quad \partial_{\Omega_-^2} \Omega^3 \quad (2.4)$$

$$q^3 \cdot n^3 = g^3 \quad \text{on} \quad \partial_N \Omega^3 \quad (2.5)$$

$$\text{tr } p^3 = 0 \quad \text{on} \quad \partial_D \Omega^3 \quad (2.6)$$

Here we denote by $\Omega^3$ the (3-dimensional) porous medium, and by $\partial_N$ and $\partial_D$ its Neumann and Dirichlet boundaries, respectively. We denote by $\partial_{\Omega_\pm^2} \Omega^3$ the boundary of $\Omega^3$ as seen from the positive (resp. negative) side of $\Omega^2$, and the outer normal vector is always denoted $n$. The Dirichlet boundary data is set to zero for notational convenience. We emphasize the structure of the governing equations as composed of a conservation law (2.1), and a constitutive (Darcy) law (2.2). In equations (2.1-2.6) we have marked variables by a superscript '3' to emphasize that they belong in 3 dimensions, the necessity of the precision will be clear below. Note that the flux from the (2-dimensional) Neumann boundary is denoted by a superscript '2'. Throughout the manuscript, we will use $\psi$ to denote right-hand sides, which with the chosen sign convention represents fluid extraction.

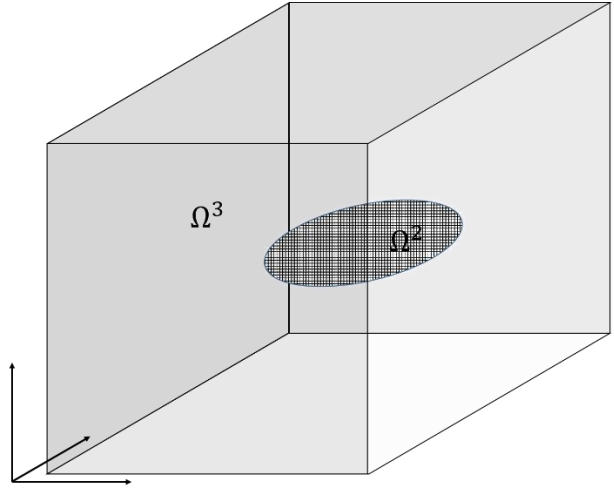

**Figure 1:** Illustration of a 3D-domain with a single 2D-fracture, see Section 4 for more examples.

Similarly, we may consider a single fracture as a (2-dimensional) manifold $\Omega^2$, whereon the governing equations can in the linearized case be expressed as [16]

$$\nabla_2 \cdot q^2 - (\lambda_+^2 + \lambda_-^2) + \psi^2 = 0 \qquad \text{on} \qquad \Omega^2 \qquad (2.7)$$

$$-\kappa_\parallel^2 \nabla_2 p^2 = q^2 \qquad \text{on} \qquad \Omega^2 \qquad (2.8)$$

$$q^2 \cdot n^2 = g^d \qquad \text{on} \qquad \partial_N \Omega^2 \qquad (2.9)$$

$$\operatorname{tr} p^2 = 0 \qquad \text{on} \qquad \partial_D \Omega^2 \qquad (2.10)$$

In equations (2.7-2.8), we denote by a double-strike the tensor operating tangentially (parallel) to the manifold and emphasize that the differential operators are 2-D by a subscript. We note that in equation (2.7), two extra terms arise. These represent the outflow from the fracture into the porous medium on the two sides of the fracture (denoted + and -). As above, fracture variables are indicated by a superscript '2' for clarity.

Considering still the case of a single fracture, equations (2.1-2.10) lead to a system of equations where $\lambda^2$ is a variable internal to the system. We thus complete the model with a constitutive law for $\lambda^2$, which takes the Darcy-like form (see e.g. [7])

$$\lambda_\pm^2 = -\kappa_\perp \big( p^2 - \operatorname{tr} p_\pm^3 \big) \qquad (2.11)$$

We remark that the within-fracture permeability $\kappa_\parallel$ and the transverse permeability $\kappa_\perp$ may in practice scale with the aperture and its inverse, respectively.

Equations (2.1-2.11) form a closed and well-posed system of equations for a porous medium including a fracture (see e.g. [8]). More generally, we note that we write these equations in a unified way, in that for $d = \{2,3\}$

$$\nabla_d \cdot q^d - \sum_{j \in \pm} \lambda_j^d + \psi^d = 0 \qquad \text{on} \qquad \Omega^d \qquad (2.12)$$

$$-\kappa_\parallel^d \nabla_d p^d = q^d \qquad \text{on} \qquad \Omega^d \qquad (2.13)$$

$$q^d \cdot n^d = \lambda^{d-1} \quad \text{on} \quad \partial_{\Omega^{d-1}}\Omega^d \quad (2.14)$$

$$\lambda_j^d = -\kappa_\perp^d\left(p^d - \operatorname{tr} p_j^{d+1}\right) \quad \text{on} \quad \partial_{\Omega^d}\Omega^{d+1} \quad (2.15)$$

$$q^d \cdot n^d = g^d \quad \text{on} \quad \partial_N \Omega^d \quad (2.16)$$

$$p^d = 0 \quad \text{on} \quad \partial_D \Omega^d \quad (2.17)$$

Equations (2.12-2.17) make sense with the convention that since there is no 4-dimensional domain in the model, the terms $\lambda^3 = 0$ and $\kappa_{\|}^3 = \kappa$.

From physical considerations, it is customary to consider all boundaries of the fracture as Neumann boundaries with $g^d = 0$, except where the boundary coincides with an outer boundary of the full domain. However, these restrictions are not necessary from a mathematical or numerical perspective, and we will retain the slightly more general formulation in order to avoid extra notation for distinguishing between internal and external boundaries of fractures.

## 2.2. Extension to general fracture configurations

Equations (2.12-2.17) are written in a way that naturally generalizes also to fracture intersections, both the 1-D line intersections as well as the 0-D point intersections of three fractures [17, 6]. We introduce some extra notation to this end. Let each domain (matrix, fracture, or intersection) be indexed by number and dimension, i.e. $\Omega_i^d$ is domain number $i \in I$, having dimensionality $d$. We consider a total of $m$ subdomains of various dimensionality. This subdivision is illustrated in Figure 2.

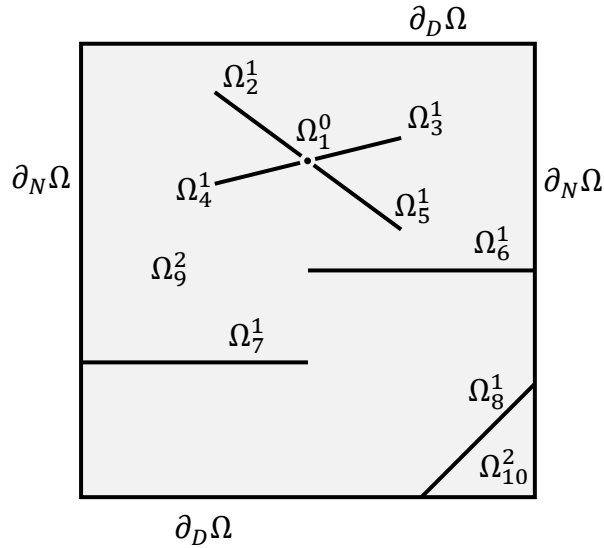

**Figure 2:** Illustration of a domain in 2D containing four fractures and an intersection, and its logical representation with two 2D-domains, seven 1D-domains, and one 0D domain.

Furthermore, let $\hat{S}_i$ be the set of neighbors of domain $i$ of dimension $d + 1$, and conversely let $\check{S}_i$ be the set of neighbors of $i$ with dimension $d - 1$. Then we can write for all $d = \{0,1,2,3\}$ and all $i \in I$ the equations

$$\nabla_d \cdot q_i^d - \sum_{j \in \hat{S}_i} \lambda_{i,j}^d + \psi_i^d = 0 \quad \text{on} \quad \Omega_i^d \quad (2.18)$$

$$-\kappa_{i,\|}^d \nabla_d p_i^d = q_i^d \quad \text{on} \quad \Omega_i^d \quad (2.19)$$

$$q_j^d \cdot n_j^d = \lambda_{i,j}^{d-1} \quad \text{on} \quad \partial_i \Omega_j^d \quad (2.20)$$

$$\lambda_{i,j}^d = -\kappa_{i,\perp}^d \left( p_i^d - \operatorname{tr} p_j^{d+1} \right) \quad \text{on} \quad \partial_i \Omega_j^{d+1} \quad (2.21)$$

$$q_i^d \cdot n_i^d = g_i^d \quad \text{on} \quad \partial_N \Omega_i^d \quad (2.22)$$

$$p_i^d = 0 \quad \text{on} \quad \partial_D \Omega_i^d \quad (2.23)$$

Note that for $d = 0$, the domain has no physical extent and no boundary such that (2.19), (2.20), (2.22), and (2.23) are void, and correspondingly $q_i^0$ is not a variable in the system. Equations (2.18-2.23) - with some variations – are equivalent or generalize the standard equations used to model fractured porous media (see [7, 1, 6] and references therein). These equations have been identified as a second-order system of mixed-dimensional partial differential equations, for which existence and uniqueness theory has been developed under fairly mild assumptions on the geometry [13]. In this work we will only consider planar fractures, but with no restrictions on their intersections or interaction with the boundary.

In order to simplify notation in the following, we consider the dimension associated with each subdomain, $d = d(i)$, to be specified, and introduce the compound variables $\mathfrak{p} = \begin{bmatrix} p_1^{d(1)} & \dots & p_m^{d(m)} \end{bmatrix}$, $\mathfrak{q} = \begin{bmatrix} q_1^{d(1)} & \dots & q_m^{d(m)} \end{bmatrix}$ and similarly for $\lambda = [\lambda_{i,j}]$. We also introduce corresponding function spaces, thus we let $\mathcal{H}_1 = \prod_i \overset{\circ}{H}_1 \left( \Omega_i^{d(i)} \right)$ and $\mathcal{L}^2 = \prod_i \prod_{j \in \hat{S}_i} L^2(\Omega_i^d(i))$. The Dirichlet boundary conditions implied by the notation $\overset{\circ}{H}_1$ only applies to the part of the boundary covered by equation (2.23).

## 2.3 Variational formulation

Before considering discretization of equations (2.18-2.23), we note that equation (2.21) is in a sense dual to the summation terms in equation (2.18), thus the system can be written as a symmetric saddle-point problem: Find $(\mathfrak{p}, \lambda) \in \mathcal{H}_1 \times \mathcal{L}^2$ such that for all $i \in I$ [from equations (2.18-2.20)]:

$$\left( \kappa_{i,\|}^d \nabla_d p_i^d, \nabla_d w_i^d \right)_{\Omega_i^d} + \sum_{j \in \hat{S}_i} \left( \lambda_{j,i}^{d-1}, \operatorname{tr} w_i^d \right)_{\partial_j \Omega_i^d} - \sum_{j \in \hat{S}_i} \left( \lambda_{i,j}^d, w_i^d \right)_{\Omega_i^d} = -(\psi_i^S, w) - \left( g_i^d, \operatorname{tr} w_i^d \right)_{\partial_N \Omega_i^d}$$

$$\text{for all} \quad w \in H_1(\Omega_i^d) \quad (2.24)$$

and [from (2.21)]:

$$\sum_{j \in \hat{S}_i} \left( p_i^d, \mu_{i,j}^d \right)_{\Omega_i} - \left( \mu_{i,j}^d, \operatorname{tr} p_j^{d+1} \right)_{\partial_i \Omega_j^{d+1}} + \left( \left( \kappa_{i,\perp}^d \right)^{-1} \lambda_{i,j}^d, \mu_{i,j}^d \right)_{\partial_i \Omega_j^{d+1}} = 0 \text{ for all} \quad \mu \in L^2(\partial_i \Omega_j^d)$$

$$(2.25)$$

By shifting indexes on the trace term in (2.24), we identify the symmetric and coupling terms as

$$a(\mathfrak{p}, \lambda; w, \mu) = \sum_{i \in I} \left( \kappa_{i,\|}^d \nabla_d p_i^d, \nabla_d w_i^d \right)_{\Omega_i^d} + \left( \left( \kappa_{i,\perp}^d \right)^{-1} \lambda_{i,j}^d, \mu_{i,j}^d \right)_{\partial_i \Omega_j^{d+1}} \quad (2.26)$$

$$b(\mathfrak{p}, \mu) = \sum_{i \in I} \left( \left( \mu_{i,j}^d, \operatorname{tr} p_j^{d+1} \right)_{\partial_i \Omega_j^{d+1}} - \sum_{j \in \hat{S}_i} \left( p_i^d, \mu_{i,j}^d \right)_{\Omega_i^d} \right) \quad (2.27)$$

For non-degenerate coefficients, equations (2.24-2.25) are well-posed by standard saddle-point theory [18], and in the remaining manuscript we will only consider this case. Nevertheless, we remark that, following similar arguments as exposed in [6], it can be shown that significant degeneracy of coefficients can be permitted, at the cost of introducing weighted spaces. In particular, it is of interest to also allow for fractures where the tangential permeability is negligible. Equations (2.24-2.25) are well-posed in this sense, since if for a given domain $\Omega_i^d$, the permeability can degenerate in the sense of $\kappa_{i,\|}^d \to 0$, as long as $\kappa_{i,\perp}^d$ remains bounded from below for all $j \in \hat{S}_i$. However, now the pressure $p_i^d$ is only in $L^2$ due to the inf-sup condition for $b(\mathfrak{p}, \mu)$ [6]. This implies that this weakly continuous formulation for fractured porous media is robust both for arbitrarily thin fractures and can also be applied to blocking fractures. We summarize the above discussion as follows:

Let an $L^2$-like norm on $\mathcal{H}_1 \times \mathcal{L}^2$ be defined as

$$\|(\mathfrak{p}, \lambda)\|^2 = \sum_{i \in I} \|p_i^d\|_{L^2(\Omega_i^d)}^2 + \sum_{j \in \hat{S}_i} \|\lambda_{i,j}^d\|_{L^2(\Omega_i^d)}^2 \tag{2.28}$$

Furthermore, let the set of indexes be refined such that $i \in I_a$ if $\kappa_{i,\|}^d > 0$ and $i \in I_b$ if $\kappa_{i,\|}^d = 0$. Then we introduce space $\mathcal{P}$ as

$$\mathcal{P} = \prod_{i \in I_a} \mathring{H}_1\left(\Omega_i^{d(i)}\right) \times \prod_{i \in I_b} L^2\left(\Omega_i^{d(i)}\right) \tag{2.29}$$

Note here that we use a circle above the function space to indicate homogeneous Dirichlet boundary conditions. Then the equations for flowing and blocking fractures can be written as find $(\mathfrak{p}, \lambda) \in \mathcal{P} \times \mathcal{L}^2$ such that

$$a(\mathfrak{p}, \lambda; w, \mu) + b(w, \lambda) - b(\mathfrak{p}, \mu) = -\sum_{i \in I}(\psi_i^s, w_i^s) - \sum_{i \in I_a}\left(g_i^d, \mathrm{tr}\, w_i^d\right)_{\partial_N \Omega_i^d} \quad \text{for all } (w, \mu) \in \mathcal{P} \times \mathcal{L}^2 \tag{2.30}$$

The solution of (2.30) is characterized by the following Lemma.

### Lemma 2.1
Equation (2.30) has a unique solution $(\mathfrak{p}, \lambda) \in \mathcal{P} \times \mathcal{L}^2$, satisfying

$$\|(\mathfrak{p}, \lambda)\| \leq C\|(\psi, g)\| \tag{2.31}$$

Provided that there exists constants $\kappa_{0,\perp}$ and $\kappa_{\infty,\perp}$ for all $i$, holds that $0 < \kappa_{0,\perp} \leq \kappa_{i,\perp}^d \leq \kappa_{\infty,\perp} < \infty$, and that

a) There is a lower bound $\kappa_{0,\|}$ such that for all $i \in I_a$, it holds that $\kappa_{i,\|}^d \geq \kappa_{0,\|} > 0$, while,
b) For all $i \in I_b$ there it holds that $j \in I_a$ for all $j \in \hat{S}_i$.

*Proof.*

For the two cases in the proof for $I_a$ and $I_b$, respectively, we indicate variables in these domains by similar subscripts. Then formally, equations (2.26) take the form

$$\begin{pmatrix} \kappa_{0,\parallel}\Delta_a & \Sigma & \Sigma & 0 \\ -\Sigma^T & \kappa_\perp^{-1} & 0 & 0 \\ -\Sigma^T & 0 & \kappa_\perp^{-1} & -\Sigma^T \\ 0 & 0 & \Sigma & 0 \end{pmatrix} \begin{pmatrix} p_a \\ \lambda_{a,a} \\ \lambda_{a,b} \\ p_b \end{pmatrix} = -\begin{pmatrix} \psi_a + g \\ 0 \\ 0 \\ \psi_b \end{pmatrix}$$

Here, $\Delta_a$ represents the $H_1$ bilinear forms on $\Omega_i^d$, $\kappa_\perp$ represents the $L^2$ bilinear forms om $\partial_j \Omega_i^d$, while $\Sigma$ are the duality pairings in (2.27). The upper-left 3x3 system is coercive due to the conditions of the proof. Furthermore, we obtain the well-posedness of the full system, since it is easy to show that the $\Sigma$ terms are inf-sup stable between $L^2$ spaces, indeed

$$\inf_{\substack{p_i^d \in L^2(\Omega_i^d) \\ i \in I_b}} \sup_{\mu_{i,j}^d \in L^2(\Omega_i^d)} \frac{\sum_{j \in \hat{S}_i} (p_i^d, \mu_{i,j}^d)_{\Omega_i^d}}{\|p_i^d\|_{L^2(\Omega_i^d)} \|\mu_{i,j}^d\|_{L^2(\Omega_i^d)}} \geq 1$$

Since one may simply choose $\mu_{i,j}^d = p_i^d$. The coercivity of the upper left 3x3 system together with inf-sup for the $\Sigma$ terms is sufficient for stability of the full system by abstract saddle-point theory [18]. □

### Remark 2.2

Lemma 2.1 is not optimal in the sense that it is fairly easy to extract $H^1$ regularities on all domains $i \in I_a$, and the restrictions on $\kappa_\perp$ can be somewhat relaxed. However, as we are primarily interested in the numerical implementation in this contribution, we have chosen to keep Lemma 2.1 as simple as possible. Readers interested in the functional analysis for equations of this type are referred to the papers referenced in the introduction.

It is important to note that the main objective of exposing the equations for flow in fractured porous media on the form (2.26-2.27), is that it highlights the specific domain-decomposition like structure of the problem. Indeed, we note that on each sub-domain (be it porous media, fracture, or fracture intersections), we have a fairly standard elliptic partial-differential equation. These are coupled via interface variables, $\lambda_{i,j}^d$. This structure is key to design general and flexible discretization approaches, as introduced in the next section.

## 3. Discretizations for fractured porous media

Our exposition of the mathematical model for fractured porous media emphasizes two main aspects of the model, namely the second-order elliptic PDE within each domain, and the flux-coupling terms. Numerous discretization methods have been constructed for second-order elliptic differential equations – many of these are bespoke to the particular challenges associated with flow in highly heterogeneous porous media (for an introduction, see the books [19, 20, 21]). Herein, we will prove that any stable discretization for flow in (fixed-dimensional) porous media can be applied to fractured porous media through the framework introduced in the preceding section.

We subdivide this section in three parts, in order to provide the mixed-dimensional flux coupling (MDFC) discretization framework, its abstract analysis, and a concrete example using finite elements.

To be precise, we consider each domain $\Omega_i^d$ and its Neumann boundary $\Gamma_i^d = \partial_N\Omega \cup_{j\in \check{S}} \partial_j \Omega_i^d$ as endowed with a numerical discretization (note that $\Gamma_i^d$ includes all boundaries to lower-dimensional manifolds). We will only consider linear discretizations, however the approach should be applicable also to non-linear discretizations (for a recent contribution in this direction from Helmig's group, see [22]). We do not require that a discrete grid be defined, however we let the discrete representation of $L^2(\Omega_i^d)$ and $L^2(\Gamma_i^d)$ be denoted as $N_h(\Omega_i^d)$ and $N_h(\Gamma_i^d)$, respectively. For domains $i \in I_a$, i.e. where the fractures are permeable with $\kappa_{i,||}^d \geq \kappa_{0,||}$, the solution operator of the numerical discretization of the heterogeneous elliptic equation on a given domain $i \in I_a$ can be stated as $\mathcal{N}_i^d : [N_h(\Omega_i^d), N_h(\Gamma_i^d)] \to [N_h(\Omega_i^d), N_h(\Gamma_i^d)]$. This solution operator maps sinks and Neumann data to pressures and pressure traces, as made precise below. Here, we recall that we for notational simplicity only consider homogeneous boundary conditions on the Dirichlet boundaries, and as such suppress the Dirichlet boundary data. For domains $i \in I_b$, the solution operator is void, as there is no differential equation on these domains.

We will use the natural requirement that the numerical discretizations provided are consistent approximations in the following sense: Let $i \in I_a$, and let $[p,t] = \mathcal{N}_i^d(\psi, \theta)$, for $(\psi, \theta) \in N_h(\Omega_i^d) \times N_h(\Gamma_i^d)$, then this quadruplet of variables approximates the solution to the elliptic differential equation

$$\nabla_d \cdot \left(-\kappa_{i,||}^d \nabla_d p\right) + \psi \approx 0 \quad \text{on} \quad \Omega_i^d \tag{3.1}$$

$$\left(-\kappa_{i,||}^d \nabla_d p\right) \cdot n - \theta \approx 0 \quad \text{on} \quad \Gamma_i^d \tag{3.2}$$

$$t - \operatorname{tr} p \approx 0 \quad \text{on} \quad \Gamma_i^d \tag{3.3}$$

$$\operatorname{tr} p \approx 0 \quad \text{on} \quad \partial_D \Omega_i^d \tag{3.4}$$

The precise interpretation of $\approx$ will depend on the chosen numerical method. We note that standard methods such as finite volume, finite element, mixed-finite element and spectral methods all fall within this framework, where the approximation implied by the $\approx$ signs of equations (3.1-3-4) can for most numerical methods be characterized by grid regularity, material parameters, grid resolution, etc. By assumption, we consider only stable numerical methods, in the sense of a negative eigenvalue-spectrum for the numerical solution operators $\mathcal{N}_i^d$, with potentially a single degenerate eigenvalue for subdomains where $\partial_D \Omega_i^d = \emptyset$, and we will denote the smallest (i.e. most negative) nondegenerate eigenvalue of $\mathcal{N}_i^d$ as $-n_i^d$. Furthermore, the system (3.1-3.4) is self-adjoint, so that in many cases the numerical method $\mathcal{N}_i^d$ will be symmetric (see Section 3.3 below for the case of finite elements).

### 3.1 MDFC: A unified discretization of fractured porous media

To provide a discretization for fractured systems, a grid $\mathcal{T}_i^d$ is introduced on the lower-dimensional manifolds $\Omega_i^d$ on which the boundary flux variables $\lambda_{i,j}^d$ will be defined. We emphasize that this mortar-like grid $\mathcal{T}_i^d$ can be chosen independently of any grid used by the numerical methods $\mathcal{N}_i^d$ and $\mathcal{N}_{\check{S}_i^d}^{d+1}$, thus we impose a minimum of restrictions on the grids. Nevertheless, note that this construction ensures that the flux variables on either side of a fracture (or either sides of fracture intersections) are conforming with each other. The precise relationships between the admissible grids $\mathcal{T}_{i,j}^d$ as implied by the numerical methods $\mathcal{N}_i^d$, will be made clear below. For the sake of symmetry, we also define grids $\mathcal{G}_i^d$ for the Neumann data on $\partial_N \Omega_i^d$.

To formulate discrete methods for fractured porous media, we represent the flux variable as piecewise constant on the mortar grid $\mathcal{T}_i^d$, thus $\lambda_{i,j}^d \in P_0(\mathcal{T}_i^d)$ and $g_i^d \in P_0(\mathcal{G}_i^d)$ (higher-order approximations are also possible, but the regularity of the problem does not seem to justify this). We introduce projection operators in order to move between the degrees of freedom of the numerical methods $\mathcal{N}_i^d$ and the mortar grids $\mathcal{T}_i^d$. We first define the compound operator projecting from the coupling variables on the mortar grids to the subdomain degrees of freedom

$$\Pi_{N_h(\Omega_i^d)}: \left[P_0(\mathcal{T}_i^d), P_0\left(\mathcal{T}_{\check{S}_i}^{d-1}\right), P_0(\mathcal{G}_i^d)\right] \to \left[N_h(\Omega_i^d), N_h(\Gamma_i^d)\right] \tag{3.5}$$

and conversely from the numerical variables to the coupling variables

$$\Pi_{P_0(\mathcal{T}_i^d)}: \left[N_h(\Omega_i^d), N_h(\Gamma_i^d)\right] \to \left[P_0(\mathcal{T}_i^d), P_0\left(\mathcal{T}_{\check{S}_i}^{d-1}\right), P_0(\mathcal{G}_i^d)\right] \tag{3.6}$$

Now, our MDFC discretization framework for fractured porous media takes the form: For given numerical discretizations $\mathcal{N}_i^d$: Find $\lambda_{i,j}^d \in P_0(\mathcal{T}_i^d)$, for all $i \in I$ and $j \in \hat{S}_i$ such that

$$\left(p_i^d, \mu_j\right)_{\partial_i \Omega_j^{d+1}} - \left(t_{i,j}^d, \mu_j\right)_{\partial_i \Omega_j^{d+1}} + \left(\kappa_{i,\perp}^{-d} \lambda_{i,j}^d, \mu_j\right)_{\partial_i \Omega_j^{d+1}} = 0 \qquad \text{for all} \quad \mu_j \in P_0(\mathcal{T}_i^d) \tag{3.7}$$

subject to the discrete constraints:

$$\left[p_i^d, t_{l,i}^{d-1}, z_i^d\right] = \Pi_{P_0(\mathcal{T}_i^d)} \mathcal{N}_i^d \left(\psi_i^d + a_i^d, b_i^d\right) \qquad \text{for all } i \in I_a \tag{3.8}$$

$$\left[a_i^d, b_i^d\right] = \Pi_{N_h(\Omega_i^d)} \left[-\sum_{j \in \hat{S}_i^d} \lambda_{i,j}^d, \lambda_{\check{S}_i,i}^{d-1}, g_i^d\right] \qquad \text{for all } i \in I_a \tag{3.9}$$

The dummy variables $a_i^d$ and $b_i^d$ have the interpretations of sinks and fluxes due to the interactions with other domains, respectively. In contrast, the variables $p_i^d$ and $t_i^d$ are the pressure and pressure traces after projection onto the grids $\mathcal{T}_i^d$. The variable $z_i^d$ is the pressure trace projected onto the Neumann boundaries, and is not used with the boundary conditions considered herein (but would be relevant with Robin boundary conditions).

This MDFC scheme has a particularly simple interpretation: For each subdomain $i \in I_a$, $\mathcal{N}_i^d$ can be interpreted as a generalized Neumann-Dirichlet map, in the sense that it maps boundary fluxes (which also take the apparent form of sources for neighboring domains of $d - 1$) to Dirichlet data (where conversely, for $d < n$, the internal values are considered Dirichlet data for neighboring domains of dimension $d + 1$). As such, equation (3.8) resolves the internal differential equations in each subdomain, equations (3.9) is the projection of variables from the flux grids to the numerical boundary (and source) data, while equation (3.7) simply states that the flux $\lambda_{i,j}$ between a fracture and its surroundings should satisfy a form of Darcy's law, depending on the difference in pressure $p$ of the fracture and the pressure $t$ at the boundary of the surroundings. Equations (3.7-3.9) are thus a Schur-complement formulation of the discrete problem.

## 3.2 Abstract analysis

Let the discretization methods corresponding to the solution operators $\mathcal{N}_i^d$ be collected in a linear system, i.e. we state equation (3.8) on the form:

$$[\mathfrak{p}, \mathfrak{t}] = \mathcal{N}[\psi + \mathfrak{a}, \mathfrak{b}] \qquad (3.10)$$

Similarly, we denote the compound projection operators $\Pi_{\mathcal{T}}$ and $\Pi_{\mathcal{N}}$. Furthermore, denote by $D$ the discrete divergence operators from equation (3.9), which sums flux variables associated with a fracture while retaining Neumann boundary data i.e.

$$D\lambda|_i = \begin{pmatrix} -\sum_{j \in \hat{S}_i^d} \lambda_{i,j}^d \\ \lambda_{i,k}^{d-1} \end{pmatrix} \qquad (3.11)$$

Finally, let the diagonal mass matrix associated with the inner product $\left(\left(\kappa_{i,\perp}^d\right)^{-1}\lambda_i^d, \mu\right)$ appearing in equation (3.9) be denoted $\kappa^{-1}$. Then we can eliminate the subdomain variables from the discrete system (3.7-3.9) to obtain a Schur-complement system only in terms of the flux variables, i.e.

$$\left(\kappa^{-1} + D^T \Pi_{\mathcal{T}} \mathcal{N} \Pi_{\mathcal{N}} D\right)\lambda = D^T \Pi_{\mathcal{T}} \mathcal{N}[\psi, \Pi_{\mathcal{N}} g] \qquad (3.12)$$

From the Schur complement form, we immediately obtain the following result:

### Lemma 3.1

Let all subdomain discretization methods $\mathcal{N}_i^d$ be negative definite for $i \in I_a$ (i.e. $\partial_D \Omega_i^d \neq \emptyset$ for all $i \in I_a$), and furthermore let the assumptions of Lemma 2.1 hold. Then if the projection operators are negative transposes, such that $\Pi_{\mathcal{T}}^T = -\Pi_{\mathcal{N}}$, the Schur-complement system (3.12) is stable, with no degenerate eigenvalues.

*Proof*: By the choice of $\lambda_{i,j} \in P_0(\mathcal{T}_i^d)$, the $\kappa^{-1}$ matrix is diagonal, and has positive eigenvalues bounded below by $\kappa_{0,\perp}^{-1}$. Thus, it is sufficient to show that the remaining term has non-negative eigenvalues. But since $\mathcal{N}$ is negative definite by the assumption of the lemma, then $D^T \Pi_{\mathcal{T}} \mathcal{N} \Pi_{\mathcal{N}} D = -(\Pi_{\mathcal{N}} D)^T \mathcal{N} \Pi_{\mathcal{N}} D$ will be non-negative definite. The result follows since the right-hand side operator is bounded by the assumption of the Lemma. □

In order to allow for fractures (and intersections, etc.) which do not have a Dirichlet boundary, the arguments of Lemma 3.1 must be refined. To this end, let $\bar{I}_a$ be the subset of $I_a$ which do not have a Dirichlet boundary. For these domains, we have a pure Neumann problem, and equations (3.8) are expected to constrain the solutions up to a constant (pressure). For the analysis, we therefore introduce an auxiliary constant pressure $\bar{p}_i^d$ for each domain $i \in \bar{I}_a$, and introduce the modified numerical methods $\widetilde{\mathcal{N}}_i^d : [N_h(\Omega_i^d), N_h(\Gamma_i^d)] \setminus \mathbb{R} \to [N_h(\Omega_i^d), N_h(\Gamma_i^d)] \setminus \mathbb{R}$, i.e., the solution corresponding to equations (3.1-3.4) with a compatibility condition (fluxes and sinks must sum to zero), and the additional constraint that the pressure has mean value zero. For $i \neq I_a \setminus \bar{I}_a$, the solution operator is unaltered, $\widetilde{\mathcal{N}}_i^d = \mathcal{N}_i^d$.

Equation (3.10) is then restated as

$$[\mathfrak{p}, \mathfrak{t}] = \widetilde{\mathcal{N}}[\psi + \mathfrak{a}, \mathfrak{b}] + \bar{\mathfrak{p}} \qquad (3.13)$$

Inserting

$$\left(\kappa^{-1} + D^{\mathrm{T}}\Pi_{\mathcal{T}}\widetilde{\mathcal{N}}\Pi_{\mathcal{N}}D\right)\lambda + \left(D^{\mathrm{T}}\Pi_{\mathcal{T}}\right)_{I,\bar{I}_a}\bar{p} = D^{T}\Pi_{\mathcal{T}}\mathcal{N}[\psi,\Pi_{\mathcal{N}}g] \tag{3.14}$$

With the compatibility constraint that

$$(\Pi_{\mathcal{N}}D)_{\bar{I}_a,I}\lambda = 0 \tag{3.15}$$

### Lemma 3.2

Let all subdomain discretization methods $\widetilde{\mathcal{N}}_i^d$ be negative definite for $i \in I_a$, and furthermore let the assumptions of Lemma 2.1 hold. Furthermore, let $I_a \setminus \bar{I}_a$ contain at least one domain. Then if the projection operators are negative transposes, such that $\Pi_{\mathcal{T}}^{\mathrm{T}} = -\Pi_{\mathcal{N}}$, the saddle-point system (3.14-3.15) is stable, with no degenerate eigenvalues.

*Proof*: By the assumptions of Lemma 3.1, $\left(D^{\mathrm{T}}\Pi_{\mathcal{T}}\right)_{I,\bar{I}_a} = -(\Pi_{\mathcal{N}}D)_{\bar{I}_a,I}$. Moreover, by similar argument to Lemma 3.1, it holds that $\left(\kappa^{-1} + D^{\mathrm{T}}\Pi_{\mathcal{T}}\widetilde{\mathcal{N}}\Pi_{\mathcal{N}}D\right)$ is coercive. It remains to show inf-sup for $\left(D^{\mathrm{T}}\Pi_{\mathcal{T}}\right)_{I,\bar{I}_a}$. I.e., we must show that

$$\inf_{\bar{p}\in\mathbb{R}^{|\bar{I}_a|}}\sup_{\lambda}\frac{(\Pi_{\mathcal{N}}D)_{\bar{I}_a,I}\lambda\cdot\bar{p}}{\|\lambda\|\|\bar{p}\|} \geq C \tag{3.16}$$

This result is obtained by considering (all) $i$ such that $i \in I_a \setminus \bar{I}_a$. Construct a rooted tree(s) $\mathfrak{T}$ from $i$ spanning all subdomains (this can always be done for connected domains). Then for leaves (i.e. terminal nodes of the tree) $j$ we set $\lambda_{j,k} = \bar{p}_j$, where $k$ is the parent of $j$ (we use the sign convention that $\lambda_{j,k} = -\lambda_{k,j}$ if $k$ is in $\check{S}_j$, and it is sufficient to consider $\lambda_{j,k}$ constant). Proceeding in this manner recursively, let $j$ be a node in the tree and let $\lambda_{j,l}$ be determined for all branches extending from $j$. Then set $\lambda_{j,k} = \bar{p}_j - \sum_l \lambda_{j,l}$. Proceeding until the root of the tree, we see by construction that $(\Pi_{\mathcal{N}}D)_{\bar{I}_a,I}\lambda = \bar{p}$, and that $\|\lambda\| \leq c\|\bar{p}\|$, where $c$ increases with the depth of the tree(s) $\mathfrak{T}$. For a finite geometry, $C$ is therefore bounded by the geometry of the fracture network, and independent of the discretization methods. The solvability and bounded eigenvalues of (3.14-3.15) then follows from standard theory [18]. □

In practice, it is of course also of interest to obtain values for the discrete solutions $p_i^d$, and not only the flux exchanges $\lambda$. This result is slightly more subtle, in a similar sense as Lemma 2.1. To prepare, we write equation (3.12) in the same form as used in the proof of Lemma 2.1.

$$\begin{pmatrix} \mathcal{A}_a & (\Pi_{\mathcal{N}}D)_{a,a} & (\Pi_{\mathcal{N}}D)_{a,b} & 0 \\ (D^T\Pi_{\mathcal{T}})_{a,a} & \kappa_a & 0 & 0 \\ (D^T\Pi_{\mathcal{T}})_{b,a} & 0 & \kappa_b & (D^T\Pi_{\mathcal{T}})_{b,b} \\ 0 & 0 & (\Pi_{\mathcal{N}}D)_{b,b} & 0 \end{pmatrix} \begin{pmatrix} p_a \\ \lambda_{a,a} \\ \lambda_{a,b} \\ p_b \end{pmatrix} = \begin{pmatrix} \phi_a \\ 0 \\ 0 \\ \phi_b \end{pmatrix} \tag{3.17}$$

Here, the linear operators $\mathcal{A}$ are the inverses of $\mathcal{N}$, and represent the linear discretizations underlying the numerical solution. Hence, equation (3.17) is also structurally similar to the natural implementation of the methodology. It is also important to note that the form (3.17) is agnostic to whether a domain is in $\bar{I}_a$, thus from the perspective of implementation, it will in many cases not be necessary to introduce special treatment of these domains as in Lemma 3.2. We now obtain a similar result as for the continuous case, in the sense that

### Theorem 3.3

Equation (3.17) is well-posed, provided that the assumptions of Lemma 2.1 and 3.1 (or 3.2) hold, and that furthermore

- c) The largest eigenvalues $n_i^d$ of the numerical methods $\widetilde{\mathcal{N}}_i^d$ are bounded from above.
- d) The discrete projection operators $\Pi_\mathcal{N}$ satisfy discrete inf-sup conditions for all pairs $i$ and $j$ appearing in condition b) of Lemma 2.1.

*Proof*: The proof is identical to Lemma 2.1 in the continuous case. □

We make the following remarks regarding Theorem 3.3 and its implications for MDFC:

1. All standard numerical methods for elliptic partial differential equations will satisfy condition c) in the theorem, thus essentially any numerical method can be applied to fractured porous media through the MDFC approach given in Section 3.1.
2. There are no restrictions on the grids $\mathcal{T}_i^d$ in relation to the numerical methods $\mathcal{N}_i^d$ as long as the fracture permeabilities $\kappa_{i,\|}^d$ do not degenerate. In particular, for grid-based numerical methods $\mathcal{N}_i^d$, non-matching grids, both coarser and finer, can be used between the external domain and $\mathcal{T}_i^d$, and furthermore into the internal domain.
3. In practice, conditions c) and d) of the theorem state that for subdomains where $\kappa_{i,\|}^d$ degenerates, the discrete representation of $p_i^d$ must not be finer than $\lambda_{i,j}^d$. This is similar to the typical conditions encountered in traditional mortar methods [15].
4. In the special case where $\mathcal{N}_i^d$ is chosen as the mixed-finite element method, analysis shows that spatially degenerating $\kappa_{i,\|}^d$ can be allowed, thus circumventing the binary structure of Lemma 2.1 and Theorem 3.3 [6].

### Corollary 3.4

A sequence of grids $\{\mathcal{T}_i^d\}_j$, numerical methods $\{\mathcal{N}_i^d\}_j$ and projection operators $\{\Pi\}_j$, where increasing $j$ is understood to enumerate finer grids, will be a convergent approximation to equation (2.30), provided the approximations to equations (2.18-2.23) are consistent.

*Proof*: Since the problem is linear, stability and consistency are sufficient for convergence. □

### 3.3 Worked example: Finite element methods

In order to make the presentation more concrete, we consider the finite element method with continuous linear Lagrange elements in the framework presented above. Thus, for each $\Omega_i^d$ let $\mathcal{U}_i^d$ be the corresponding grid, with nodal degrees of freedom.

Then for $i \in I_a$, the elements of the sub-matrices $A_i^d$ of $\mathcal{A}$ are simply given by the inner products of $p, w \in P_1(\mathcal{U}_i^d)$

$$\left(\kappa_{i,\|}^d \nabla_d p_i^d, \nabla_d w\right)_{\Omega_i^d} \tag{3.17}$$

with Neumann data implemented as natural boundary conditions through the duality pairing

$$\sum_{j \in \check{S}_i} \langle \lambda_{j,i}^{d-1}, \mathrm{tr}\, w \rangle_{\partial_j \Omega_i^d} \tag{3.18}$$

The Neumann boundary conditions are exactly dual to the evaluation of traces, and thus the operator $\mathcal{N}_i^d$ will be self-adjoint. Standard finite element theory further guarantees that the required bound on the eigenvalues of $\mathcal{N}_i^d$ holds independent of grid spacing with [23]

$$n_i^d \leq C \left(\kappa_{i,\|}^d\right)^{-1} \tag{3.19}$$

Since the solution $p_i^d$ and its trace live in finite-dimensional subspaces of $L^2$, the projection operators become defined in the standard way, i.e. for $\lambda_{i,j} \in P_0(\mathcal{T}_i^d)$ the projection $\Pi_{\mathcal{N}_i^d} \lambda_{i,j} \in P_1(\mathcal{U}_i^d)$ satisfies

$$\left(\Pi_{\mathcal{N}_i^d} \lambda_{i,j}, v\right) = (\lambda_{i,j}, v) \qquad \text{for all} \quad v \in P_1(\mathcal{U}_i^d) \tag{3.20}$$

It is therefore clear that $\Pi_{\mathcal{T}}^T = \Pi_{\mathcal{N}}$. Thus, all the conditions of Theorem 3.3 are satisfied, provided that the grids $\mathcal{U}_i^d$ are no finer than $\mathcal{T}_i^d$ whenever $\kappa_{i,\|}^d \to 0$.

We note that the finite element approximation could also be obtained directly from Section 2 by simply using the finite-dimensional spaces and the bilinear forms defined in equations (2.26-2.27). Thus equations (3.7-3.9) with the numerical methods $\mathcal{N}_i^d$ defined by equations (3.14-3.15) and projection operators defined by equation (3.15) is equivalent to the symmetric and bilinear saddle-point problem:

Find $\left(p_{i,h}^d, \lambda_{i,h}^d\right) \in P_1(\mathcal{U}_i^d) \times P_0(\mathcal{T}_i^d)$ such that

$$a\left(p_{i,h}^d, \lambda_{i,h}^d, w, \mu\right) + b'\left(p_{i,h}^d, \mu\right) + b'\left(q, \lambda_{i,h}^d\right) = (\psi_i^s, w) \quad \text{for all } (w, \mu) \in P_1(\mathcal{U}_i^d) \times P_0(\mathcal{T}_i^d) \tag{3.21}$$

This discretization is consistent within each domain (for shape-regular grids), thus it represents a consistent approximation to equations (2.30) whenever the boundary data is resolved. Note that for matching grids between the mortar space and the finite element spaces, this does *not* hold, since a "checkerboard"-type oscillation in $L^2$ is projected to zero by $\Pi_{\mathcal{N}_i^d}$. Thus while the lowest-order finite element variant of MDFC is stable for matching grids, it requires that the grids $\mathcal{T}_i^d$ are coarser than the grids chosen for resolving the elliptic partial differential equations in order to be a convergent numerical discretization for equations (2.18-2.23).

While the approach as stated above is sufficient, in the sense of obtaining a stable and convergent discretization, we also remark that an improved method would likely be obtained by honoring the structure of $\mathcal{P}$ from section 2.3, and thus using $p_{i,h}^d \in P_0(\mathcal{T}_i^d)$ for $i \in I_b$. In particular, this would eliminate the projection errors associated with the low-permeable fractures. This highlights the flexibility of the framework to accommodate different discretizations in the different domains, bespoke to the physical processes.

## 4. Example calculations

To confirm the theory derived above, we propose two synthetic test cases in which the ambient space is two- and three-dimensional, respectively. Out of the range of numerical methods to which the MDFC applies, we consider five discretization schemes, summarized below.

Two mixed methods are employed, namely the mixed finite element (RT0), and the dual virtual element method (VEM). The mixed finite element, considered in [6], is given by Raviart-Thomas elements of lowest order for the fluxes and piecewise constants for the pressure in all dimensions. On the other hand, VEM [9] employs a single degree of freedom per face for the fluxes without explicitly specifying the basis functions and represents pressures as piecewise constants. Thirdly, employing nodal-based, linear Lagrange elements in all dimensions leads to the primal formulation (P1) as presented in Section 3.3. This is the only method considered in this work which does not respect local mass conservation. Finally, two finite volume methods are considered, the two-point flux approximation (TPFA) and the multi-point flux approximation scheme (MPFA) [24].

In line with the spirit of the theory presented in this work, the coupling between dimensions employs a flux mortar variable, defined as piecewise constants on a separately generated, lower-dimensional grid. All computations are performed using the open-source simulation tool PorePy [25, 26].

### 4.1 Two-dimensional fracture system

The first example, obtained from [6], consists of a unit square with five one-dimensional fractures as given in figure 2. Immersed in the top half of the domain are two intersecting, conductive fractures with permeability $\kappa_\perp = 10^4$ and $\kappa_{||} = 1$. Note that due to the dimensionless scaling, this corresponds to fractures that are equally conductive in the parallel direction (in terms of volume per unit pressure drop) to the full porous unit square domain. Below are two half-immersed blocking fractures ($\kappa_\perp = 1, \kappa_{||} = 10^{-4}$) and finally, a conductive fracture separates the lower right corner. The boundary conditions are chosen as a unit pressure drop from top to bottom and no-flow conditions on the sides. The matrix permeability is set to 1.

This example is designed to contain all the elements that constitute challenges for numerical methods for fractured porous media: The two intersecting fractures represent both 1D and 0D domains which have no contact with the boundary, thus the numerical methods $\mathcal{N}_i^d$ on these domains will contain a degenerate eigenvalue (i.e. the pressure solutions are only defined up to a constant). Moreover, the low-permeable and horizontal fractures are expected to lead to singularities in the solution in the 2D domain. Finally, in the lower corner there is a domain which intersects both a Dirichlet and a Neumann boundary.

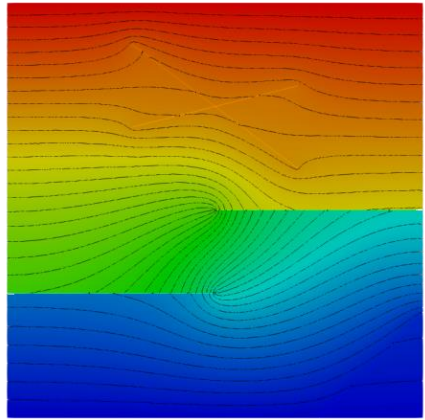

**Figure 3:** The contour lines and color scale of the reference solution on the domain given in Figure 2. The different qualitative aspects of the solution between the conductive and blocking fractures can be clearly seen.

In terms of mesh generation, the one-dimensional fracture grids match the trace of the adjacent two-dimensional grids. The mortar grid is then constructed at each fracture to have approximately 75% of the number of elements compared to the inner, lower-dimensional mesh.

Qualitatively, all numerical methods produce the same pressure distributions. Aside from artifacts due to the coarseness of the grid, all methods produce solutions which are visually indistinguishable from the figure 3. We turn to a more quantitative measure in order to expose differences between the discretizations. Since the only common property between the methods is the mortar variable, we compute its $L^2$-error with respect to a fine-scale solution obtained using the RT0 method. In case of convergence, the rate will be limited to first order with respect to the mesh size, since the mortar variable is represented by piecewise constants.

The results of this convergence test are shown in figure 4. For the one-dimensional mortar variables, very similar behavior is observed for the methods RT0, VEM, and MPFA, exhibiting stable and linear convergence. The two remaining methods show lower than first-order convergence on average. For P1, we speculate that this is due to its lack of local mass conservation, since the error is measured in a flux variable. For TPFA, this deviation is likely due to the lack of consistency in the method (i.e. the approximation error to equations (3.1-3.4) does not necessarily go to zero with grid size). We emphasize that all methods are robust and stable from a linear algebra perspective on all grids.

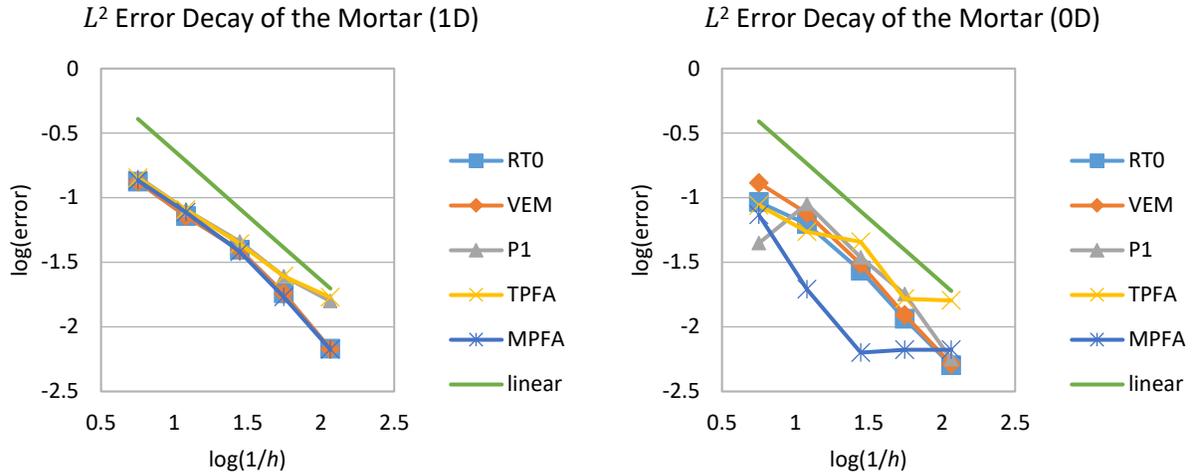

**Figure 4**: The $L^2$ errors in the mortar variable decrease with the mesh size for this range of $h$.

The error in the mortar variables defined at the zero-dimensional intersection is analyzed in figure 4 (right). These results are slightly more sporadic since an accumulation of errors can occur from the higher dimensions, and since this essentially represents a point evaluation of the solution. Moreover, the grids used in the computations are not nested and mesh sensitivities of the method with may be the cause of these effects. Nevertheless, the overall trend in all methods is a decrease in error as the mesh becomes finer. The mixed finite element methods exhibit a more monotone decay in comparison to the finite volume methods, likely due to the fact that the reference solution is calculated with the RT0 method on a finer grid.

### 4.2 Stability

It is of interest to verify the claims of Theorem 3.3. In particular, we wish to address whether the discrete representation leads to a linear system which has a lower bound on condition numbers, which is independent of grid resolutions for non-degenerate parameters, and allows for degenerate parameters in the sense of conditions a)-d) in the proof. We have chosen the condition number of the Schur complement system (3.12) as a proxy for the stability of the method, arguing (as in the preceding section) that the condition number of the full system will depend on the particular features of the numerical methods and grids utilized outside of the fractures to an extent where it is difficult to make a fair comparison.

In order to emphasize grids and parameters, we simplify the example from Section 4.1 by omitting the fractures which do not touch the boundary, and replacing the no-flow boundary conditions on the sides of the domain by a linear pressure variation. We can then consider Theorem 3.1 purely in terms of the mortar variables $\lambda_{i,j}$. Furthermore, in order to reduce the parameter space, we will let the remaining three fractures have the same parameters $\kappa_\perp$ and $\kappa_{||}$.

We fix the grid in the 2D domain with a resolution corresponding to the second-coarsest grid (approximately 4.5k triangles) in the convergence test of Section 4.1. Then in addition to the two fracture parameters, we introduce two grid parameters: The relative resolution of the outer grid to the mortar grid, and the relative resolution from the mortar grid to the fracture grid.

Our aim is to see how the lowest eigenvalue of the discrete Schur-complement system (3.12) depends on the fracture parameters and grid parameters. To this end, we have conducted a suite of simulations for all methods, exploring the full 4D parameter space. We observe that the results are completely independent of $\kappa_{||}$ and the ratio of the mortar grid to the inner grid. When varying the perpendicular permeability $\kappa_\perp$, the results depend primarily on whether the mortar grid is finer or coarser than the outer grid, and weakly depends on the ratio. These results are summarized in Figure 5.

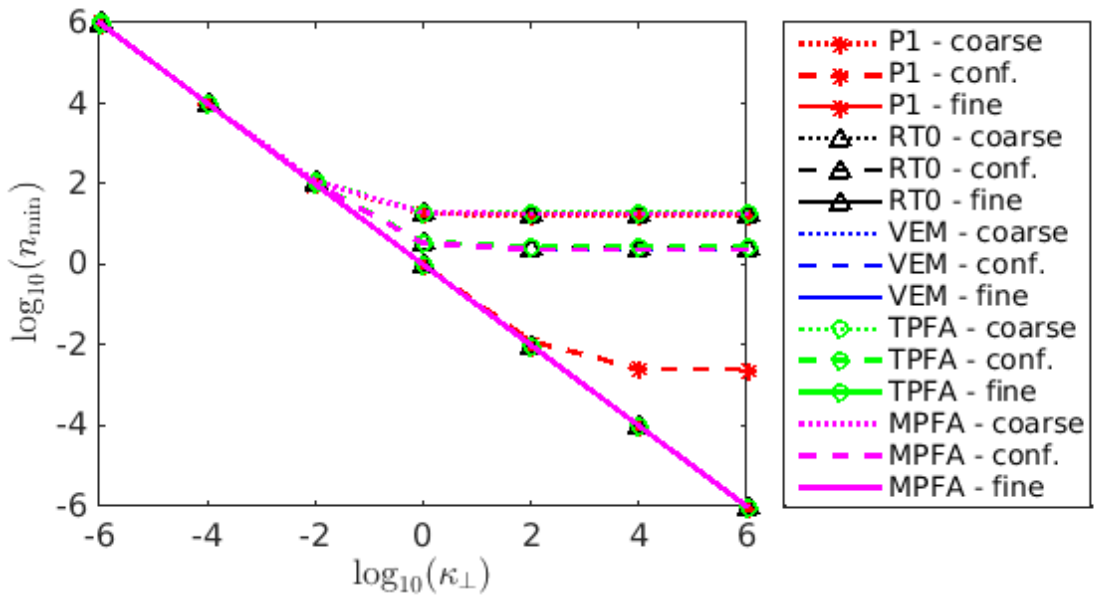

**Figure 5:** The minimum eigenvalue $n_{\min}$ of the Schur complement equation (3.12) is plotted against $\kappa_\perp$ for all five ambient numerical methods, in the case of both a finer, conforming, and coarser mortar grid (with respect to the outer grid). In all cases, the results are independent of the coarsening/refinement ratio.

From Theorem 3.3, the expected results are that the minimum eigenvalue should scale linearly with $\kappa_\perp^{-1}$. Indeed, this is what is observed for all methods in the case of small values of $\kappa_\perp$. Moreover, all methods are also stable for coarse mortar grids for large values of $\kappa_\perp$. This result reflects the fact that for coarse mortar grids, the Neumann-Dirichlet maps stabilize the system, and that numerically there is an inf-sup condition on $\Pi_\mathcal{N} D$ such that $D^T \Pi_\mathcal{T} \mathcal{N} \Pi_\mathcal{N} D$ has a lowest eigenvalue. We note however, that this does not hold for the continuous system given in equation (2.30), since the trace spaces for the pressure are not rich enough to control the mortar space. This explains why stability is lost on fine mortar grids for all methods, and is also reflected for the P1 variant of the method, which has a worse stability constant for high tangential permeability even for matching grids (see also discussion in section 3.3).. Thus, in all cases and for all grids, the MDFC method is stable, with eigenvalue bounded from below by the continuous problem.

Based on these computations, we summarize that for non-degenerate parameters, all discretizations lead to stable systems for the mortar variable, independent of grid resolution between matrix, flux-variable, and the fractures. For degenerate fracture flow $\kappa_{||}$, all methods remain stable. Finally, for degenerate fracture cross-flow $\kappa_\perp$, the results are in accordance with Theorem 3.3.

## 4.3 Three-dimensional Example

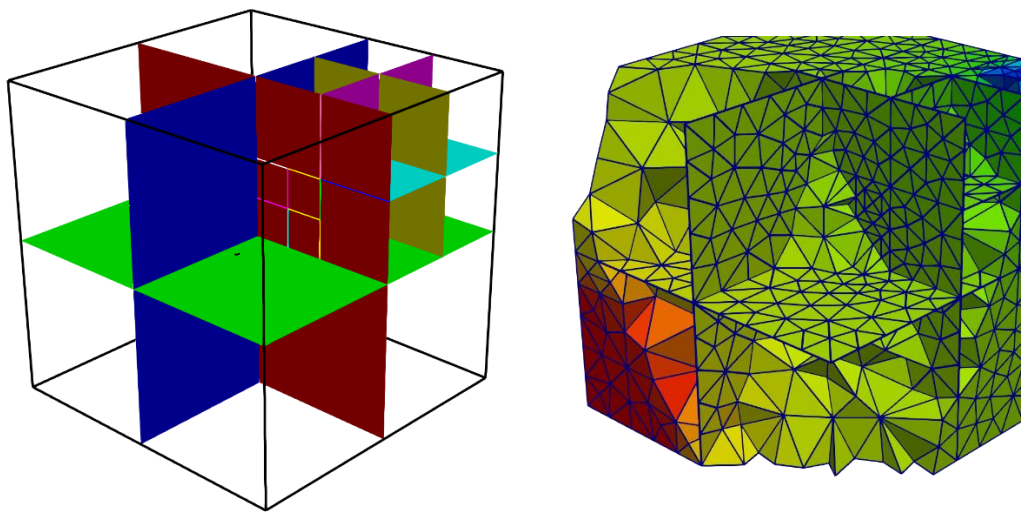

Figure 6: On the left a) the geometry of the example. On the right b) the pressure computed with RT0.

Finally, we consider simulations in a 3D problem. The computational domain is taken as the unit cube, and the fracture network for this example is reported in Figure 6 (left). The latter consist of 9 fractures with a structure similar to the Benchmark 1 in [Flemisch2017], extended to 3D. The matrix permeability is the identity tensor. We introduce the scaling factor $\zeta = 10^{-4(3-d)}$, for each lower dimensional object the normal permeability is given by $\kappa_\perp = 10^4/\zeta$ and the tangential by $\kappa_{||} = 10^4 \zeta$. Flow is forced diagonally across the domain by specifying a pressure value of 1 at boundaries characterized by $(x, y, z) < 0.4$, and similarly a pressure of $-1$ at boundaries with $(x, y, z) > 0.8$. On all other boundaries, no-flow conditions are assigned. For illustration, the numerical solution computed using RT0 is reported in Figure 5 (right).

To compare the numerical schemes, we investigate numerical convergence of the mortar variables in the same way as in Section 4.1. Three simplex grids are considered, with cell counts of about 3.5k, 4.5k and 10k tetrahedrals, together with a suitable number of triangles, line elements and points. For simplicity, we consider only matching grids in this case. Since P1 is not convergent for matching grids (see Section 3.3. and the discussion in Section 4.2), we exclude this variant of MDFC from our results. Errors in the mortar variables are computed relative to a reference solution obtained with RT0 on a grid with about 190k tetrahedral cells. The resulting error decay is depicted in Figure 7. The simulation confirms the findings in section 4.1: MPFA, RT0 and VEM all exhibit at least first order convergence for all dimensions, while TPFA again suffers from lack of consistency on the ambient grid, thus the low accuracy of the numerical method pollutes the flux variable.

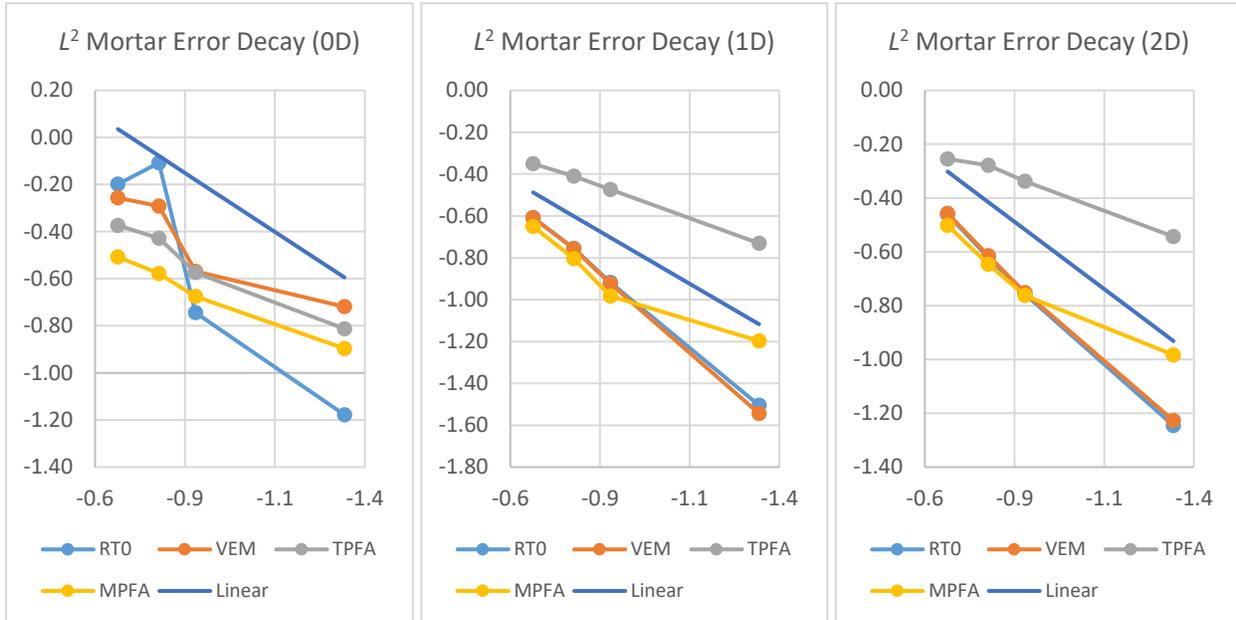

Figure 7: $L^2$ Error decay in the mortar variable for the 3d simulation reported in section 4.3.

## 5. Conclusions

We have developed a new, unified, approach to discretizing fractured porous media, termed Mixed-Dimensional Flux Coupling. The MDFC approach allows for arbitrary numerical discretizations to be used both for the porous media and the fractures. We have supported the development by both theoretical analysis, as well as numerical examples using five different numerical methods.

Several of the limitations included in this work appear to be possible to overcome. In particular, we expect that extension to non-linear discretizations [22] to be straight-forward in practice. Moreover, due to being agnostic of the numerical methods used, our theoretical results are not optimal nor exhaustive, and a more explicit treatment of the precise characteristics of the numerical methods chosen for the various components of the problem is known to provide more nuanced results [6].

In applications, coupled problems are of particular interest. In particular, the fluid flow is often coupled to transport of either mass or energy. Preliminary work in this direction is ongoing, and we expect that the MDFC framework proposed herein will accommodate such coupled problems.

We conclude by noting the importance of open-source code availability. The methods developed herein have been implemented in PorePy, and both methods and the scripts used to generate the presented results are available in the public domain at time of publication [26].

# Acknowledgements

This work has been funded in part by Norwegian Research Council grant 250223 and 244129/E20 (through the ENERGIX program).